\documentclass[11pt]{gtart}
\usepackage{amsmath, amssymb, graphicx, graphs, psfrag}

\newtheorem{thm}{Theorem}[section]

\newtheorem{cor}[thm]{Corollary}
\newtheorem{lem}[thm]{Lemma}
\newtheorem{prop}[thm]{Proposition}

\theoremstyle{definition}
\newtheorem{defn}[thm]{Definition}
\newtheorem{conj}[thm]{Conjecture}

\newtheorem{rem}[thm]{Remark}

\numberwithin{equation}{section}

\newcommand{\al}{\alpha}

\newcommand{\ga}{\gamma}
\newcommand{\Ga}{\Gamma}

\newcommand{\de}{\delta}
\newcommand{\ep}{\epsilon}

\newcommand{\Si}{\Sigma}

\newcommand{\ze}{\zeta}
\newcommand{\bfz}{{\mathbb {Z}}}

\newcommand{\x}{\times}
\newcommand{\s}{\mathbf s}

\renewcommand{\t}{\mathbf t}
\newcommand{\T}{\mathbb T}

\newcommand{\FF}{\mathcal F}

\newcommand{\Z}{\mathbb Z}

\newcommand{\R}{\mathbb R}

\newcommand{\del}{\partial}
\newcommand{\hra}{\hookrightarrow}

\newcommand{\hf}{{{\widehat {HF}}}}

\DeclareMathOperator{\tb}{tb}

\DeclareMathOperator{\id}{id}
\DeclareMathOperator{\tr}{tr}

\DeclareMathOperator{\Tor}{Tor}
\DeclareMathOperator{\rank}{rk}

\DeclareMathOperator{\Spin}{Spin}

\begin{document}

\title{Contact Ozsv\'ath--Szab\'o Invariants and\\ Giroux Torsion}

\author{Paolo Lisca}
\address{Dipartimento di Matematica ``L. Tonelli''\\
Universit\`a di Pisa\\
Largo Bruno Pontecorvo, 5\\
I-56127 Pisa, ITALY} 

\author{Andr\'{a}s I. Stipsicz}
\address{R\'enyi Institute of Mathematics\\
Hungarian Academy of Sciences\\
H-1053 Budapest\\ 
Re\'altanoda utca 13--15, Hungary}

\begin{abstract}
In this paper we prove a vanishing theorem for the contact Ozsv\'ath--Szab\'o invariants of certain contact 3--manifolds having positive Giroux torsion. We use this result to establish similar vanishing results for contact structures with underlying 3--manifolds admitting either a torus fibration over $S^1$ or a Seifert fibration over an orientable base. We also show -- using standard techniques from contact topology -- that if a contact 3--manifold $(Y,\xi)$ has positive Giroux torsion then there exists a Stein cobordism from $(Y,\xi)$ to a contact 3--manifold $(Y,\xi')$ such that $(Y,\xi)$ is obtained from $(Y,\xi')$ by a Lutz modification.
\end{abstract}


\maketitle

\section{Introduction}\label{s:intro}

In \cite{Gtor} Giroux introduced the important invariant $\Tor (Y,
\xi )$ of a contact 3--manifold $(Y, \xi )$, which is now called the
\emph{Giroux torsion}, and is defined as follows:
$\Tor (Y, \xi )$ is the supremum of the integers $n\geq 1$ for which
there is a contact embedding of 
\[
\T_n := (T^2\times [0,1],
\ker (\cos (2\pi nz)dx - \sin (2\pi nz)dy))
\]
into $(Y, \xi )$. We say that $\Tor (Y, \xi )=0$ if no such embedding
exists. 

Closed, toroidal 3--manifolds carry infinitely many universally tight contact
structures obtained by inserting copies of $\T_n$ around incompressible tori~\cite{Colin1, Colin2, CGH, HKM}. Remarkably, as the following result shows, embedded copies of $\T_n$ are the only source of infinite families of distinct tight contact structures on a closed 3--manifold.

\begin{thm}[\cite{CGH}, Theorem~1.4]\label{t:tors-quoted}
Let $Y$ be a closed 3--manifold. For every natural number $n$ the
3--manifold $Y$ carries at most finitely many isomorphism classes of
tight contact structures with Giroux torsion bounded above by $n$. \qed
\end{thm}

Contact Ozsv\'ath--Szab\'o invariants are very useful tools in studying
contact structures on closed 3--manifolds. The nature of these invariants
is still unclear though, and it is natural to ask how the invariants change under the introduction of Giroux torsion. Based on the analogy between Seiberg--Witten and Heegaard Floer theories the following is expected:

\begin{conj}\label{c:torsion}
If $\Tor (Y,\xi)>0$ then the untwisted contact Ozsv\'ath--Szab\'o
invariant $c(Y,\xi)$ vanishes.
\end{conj}

Conjecture~\ref{c:torsion} has been verified by Paolo Ghiggini~\cite{Gh3} for a class of contact Seifert fibered 3--manifolds. In the following we will extend his result and verify Conjecture~\ref{c:torsion} for a much wider family of contact 3--manifolds.

A recent result of D. Gay~\cite{Gay} asserts that a contact structure with positive
Giroux torsion is not strongly fillable. Since strongly fillable
contact structures have nonvanishing contact Ozsv\'ath--Szab\'o
invariants, Conjecture~\ref{c:torsion} is consistent with Gay's result. 
Indeed, starting from any contact 3--manifold with positive Giroux torsion, Gay constructs a symplectic cobordism which contains homologically essential 2--spheres with self--intersection zero. This suggests that such a cobordism could be used to 
find constraints on certain Seiberg--Witten invariants, thus re--proving Gay's nonfillability result. On the other hand, at present we do not understand well enough how the contact invariants of Ozsv\'ath and Szab\'o's behave under the maps induced between the relevant Heegaard Floer groups by general symplectic cobordisms. This is the main obstacle which prevents us from proving Conjecture~\ref{c:torsion} using Gay's construction. 

In this paper we build a different type of cobordism on a contact 3--manifold with positive Giroux torsion. Our cobordism is better suited than the one of Gay's to study 
the contact Ozsv\'ath--Szab\'o invariants because it is a union of Stein 2--handles, and the behaviour of the invariants under the corresponding maps is well understood. It follows that the contact invariant is always in the image of such a map, which is very useful. In fact, in the cases considered in this paper we prove that the invariant is equal to zero by showing that a certain map induced by the cobordism vanishes.

Throughout the paper, every 3--manifold will be considered to be oriented
and every contact structure positive. Recall that in \cite{OSzF1} a
variety of homology groups -- the Ozsv\'ath--Szab\'o homologies --
and a natural map
\[
\varphi_{(Y,\t)} \colon HF ^{\infty}(Y, \t )\to HF^+ (Y, \t)
\]
are defined for closed, oriented spin$^c$ 3--manifolds.  (For
more about Ozsv\'ath--Szab\'o homologies see Section~\ref{s:second}.)
\begin{defn}
We say that the closed 3--manifold $Y$ has \emph{simple
Ozsv\'ath--Szab\'o homology} at the spin$^c$ structure $\t\in
Spin^c(Y)$ if the map $\varphi_{(Y,\t)}$ is surjective.  $Y$ is called
~\emph{OSz--simple} if $Y$ has simple Ozsv\'ath--Szab\'o homology for
every spin$^c$ structure $\t\in Spin^c(Y)$.
\end{defn}

A rational homology sphere is called an \emph{$L$--space} in
\cite{OSzlens} provided that the map $\varphi_{(Y,\t)}$ is surjective
for every $\t$.  Examples of $L$--spaces can be produced by
considering plumbings of spheres along trees with no ``bad
vertices''~\cite{OSzplum} or by taking double branched covers of $S^3$
along nonsplit, alternating links~\cite[\S 3]{OSzbran}. If one uses
$\Z/2\Z$--coefficients then Seifert fibered 3--manifolds over an
orientable base with sufficiently large background Chern numbers are
OSz--simple, cf. Section~\ref{s:second}.

Given a contact 3--manifold $(Y,\xi)$, we shall denote by $\t_\xi$ the
spin$^c$ structure induced on $Y$ by $\xi$. Our first result is:

\begin{thm}\label{t:main}
Let $(Y,\xi)$ be a contact 3--manifold such that $Y$ is OSz--simple. If 
$\Tor(Y, \xi) > 1$ then $c(Y,\xi)=0$. If $b_1(Y)\leq 1$ then $\Tor(Y, \xi )>0$ already implies $c(Y, \xi )=0$.
\end{thm}

\begin{rem}\label{r:sharp}
The proof of Theorem~\ref{t:main} works under the weaker assumption
that $Y$ has simple Ozsv\'ath--Szab\'o homology at the spin$^c$
structure $\t_\xi$.  One way to check that $Y$ has simple
Ozsv\'ath--Szab\'o homology at $\t$ is to prove that $Y$ is
OSz--simple, cf.~Proposition~\ref{p:ell-par}.
\end{rem}

The following two results deal with many cases where the underlying manifolds are
not OSz--simple.

\begin{thm}\label{t:bundle}
Let $Y$ be a closed 3--manifold which admits a torus fibration over
$S^1$. Then there exists an integer $n_Y\geq 0$ such that 
for every contact structure $\xi$ on $Y$  
with 
$\Tor(Y, \xi ) > n_Y$ we have $ c(Y, \xi )=0$. 
If the monodromy $A$ of the torus fibration is trivial then $n_Y=0$.
If $A$ is elliptic ($|\tr(A)|<2$) or parabolic ($|\tr(A)|=2$)
then $n_Y\leq 1$.
\end{thm}

\begin{thm}\label{t:sfs}
Let $(Y,\xi)$ be a contact 3--manifold such that $Y$ admits a Seifert
fibration over an orientable base. If $\Tor(Y, \xi) > 2$ then using
$\Z/2\Z$--coefficients the contact Ozsv\'ath--Szab\'o invariant
$c(Y,\xi )$  vanishes.
\end{thm} 

In a slightly different direction, these vanishing results can be used to
study strong and Stein fillability of contact structures.  In~\cite{Gh2}
Ghiggini found the first examples of strongly fillable contact structures
which are not Stein fillable. His examples live on the Brieskorn 3--spheres
$-\Si(2,3,12n+5)$, $n\geq 1$. The following consequence of
Theorem~\ref{t:main}, pointed out to us by Paolo Ghiggini, slightly
generalizes~\cite[Theorem~1.5]{Gh2}.

\begin{thm}\label{t:ghiggini+}
Let $\Si$ be the Brieskorn homology 3--sphere of type $(2,3,6n+5)$,
oriented as the link of the corresponding isolated singularity. For
every $n\geq 2$, the oriented 3--manifold $-\Sigma$ carries a strongly
fillable contact structure which is not Stein fillable.
\end{thm}

Finally, using standard techniques from contact topology we establish
the following Theorem~\ref{t:reducetorsion}, which lies within the
circle of ideas of this paper and appears to be of independent
interest.  It is worth pointing out that we did not use
Theorem~\ref{t:reducetorsion} to prove any of the previous results,
except for the second part of the statement of Theorem~\ref{t:main}.
Given a smoothly embedded torus $T\subset (Y,\xi)$ with characteristic
foliation made of simple closed curves, we shall call the insertion of
a copy of $\T_1$ around $T$ a~\emph{Lutz modification} of $\xi$ along
$T$.

\begin{thm}\label{t:reducetorsion}
Let $n\geq 1$, and suppose that $\T_n$ embeds inside the contact
3--manifold $(Y,\xi)$. Then there is a sequence of Legendrian
surgeries on $(Y,\xi)$ which yields a contact 3--manifold $(Y,\xi')$
such that $(Y,\xi)$ is obtained from $(Y,\xi')$ by a Lutz modification
along an embedded copy of $\T_{n-1}$~\footnote{When $n=1$, this is to be
interpreted as an embedded 2--torus with characteristic
foliation made of simple closed curves.}.
\end{thm}

The paper is organized as follows. Section~\ref{s:second} is devoted
to the recollection of basic facts regarding Ozsv\'ath--Szab\'o
homologies and contact Ozsv\'ath--Szab\'o invariants.  We also compute
the Ozsv\'ath--Szab\'o homology groups of some of the 3--manifolds
which will appear in later arguments.  In Section~\ref{s:third} we
prove a few auxiliary results which will be used in the proofs of the
results stated above. In Section~\ref{s:fourth} we prove all the
results except Theorem~\ref{t:reducetorsion}, which is proved in
Section~\ref{s:appendix}.

{\bf {Acknowledgements}}: The second author was partially 
supported  by OTKA T49449. The authors also acknowledge partial 
support by the EU Marie Curie TOK program BudAlgGeo.
We would like to thank Paolo Ghiggini for helpful discussions.

\section{Contact Ozsv\'ath--Szab\'o invariants}
\label{s:second}

\sh{Ozsv\'ath--Szab\'o homologies} 

In the seminal papers \cite{OSzF1, OSzF2} a collection of homology groups -- the Ozsv\'ath--Szab\'o homologies -- $\hf (Y, \t), HF^{\pm}(Y,t)$ and $HF ^{\infty }(Y,\t)$ have been assigned to any closed, oriented spin$^c$ 3--manifold $(Y,\t)$.  A spin$^c$ cobordism $(X, \s )$ from $(Y_1, \t _1)$ to $(Y_2,\t_2)$ induces  $\Z [U]$--equivariant homomorphisms ${\widehat{F}}_{W, \s}$, $F^{\pm }_{W, \s }$ and $F^{\infty }_{W, \s }$ between the corresponding groups. For a fixed spin$^c$ structure $\t \in Spin^c (Y)$ these groups fit into long exact sequences
\[
\ldots \to HF^-_d (Y, \t )\to HF ^{\infty }_d (Y, \t )
\stackrel{\varphi _{(Y, \t )}}{\longrightarrow} HF^+_d (Y, \t ) \to \ldots
\]
\[
\ldots \to \hf _d (Y, \t )\to HF ^{+ }_d (Y, \t )
\stackrel{\cdot U}{\longrightarrow} HF^+ _{d-2} (Y, \t ) \to \ldots
\]
These exact sequences are functorial with respect to the maps induced
by spin$^c$ cobordisms. Throughout the paper we shall use
Ozsv\'ath--Szab\'o homology groups with $\Z$--coefficients, with the
exceptions of Theorem~\ref{t:sfs} and Proposition~\ref{p:seif},
where $\Z/2\Z$--coefficients are applied.

Ozsv\'ath--Szab\'o homology groups and the maps induced by the
cobordisms form a TQFT in the sense that the composition of two
spin$^c$ cobordisms $(W_1, \s_1)$ and $(W_2, \s _2 )$ induce a map
which can be given by the composition of the maps. There is, however,
a subtlety following from the fact that the spin$^c$ structures $\s
_i$ on $W_i$ ($i=1,2$) do not uniquely determine a spin$^c$ structure
on the union $W_1\cup W_2$. Consequently the composition formula reads
as follows:

\begin{thm}[\cite{OSzF4}, Theorem~3.4]\label{t:composition}
Suppose that $(W_1, \s_1)$ and $(W_2, \s _2 )$ are spin$^c$ cobordisms
with $\del W_1=-Y_1\cup Y_2$, $\del W_2=-Y_2\cup Y_3$ and set
$W=W_1\cup_{Y_2} W_2$. Let ${\mathcal {S}}$ denote the set of spin$^c$
structures on $W$ which restrict to $W_i$ as $\s_i$ for $i=1,2$. Then
\[
F_{W_2,\s _2}\circ F_{W_1, \s _1}=\sum _{\s \in {\mathcal {S}}} \pm
F_{W, \s }.
\]
\qed
\end{thm}

We shall say that a closed 3--manifold $Y$ has \emph{simple
Ozsv\'ath--Szab\'o homology} in the spin$^c$ structure $\t$ if the map
\[
\varphi_{(Y,\t)}\co HF^{\infty }(Y,\t) \to H^+(Y,\t)
\]
is onto.  This condition is equivalent to the requirement that
$HF_{red} (Y, \t )=0$.  We shall say that the 3--manifold $Y$ is
\emph{OSz--simple} if $Y$ has simple Ozsv\'ath--Szab\'o homology for
every spin$^c$ structures $\t\in\Spin^c(Y)$.  Since $HF_{red} (Y, \t
)\cong HF _{red} (-Y, \t)$, it follows that a closed, oriented
3--manifold $Y$ is OSz--simple if and only if $-Y$ is OSz--simple,
cf. \cite{OSzF4}.  An OSz--simple rational homology sphere is called
an \emph{$L$--space} in \cite{OSzlens}.

An important ingredient in our subsequent discussions is

\begin{prop}[\cite{OSzF4}, Lemma~8.2]\label{p:cobordism}
Suppose that $(W,\s)$ is a 4--dimensional spin$^c$ cobordism between
$(Y_1, \t _1)$ and $(Y_2, \t _2)$. If $b_2^+(W)>0$ then the map
\[
F^{\infty }_{W, \s}\colon 
HF^{\infty }(Y_1, \t _1)\to HF^{\infty } (Y_2, \t_2)
\]
is zero. \qed
\end{prop}

\begin{cor}\label{c:vanish}
Suppose that the 4--dimensional cobordism $W$ between $Y_1$ and $Y_2$
has $b_2^+(W)>0$ and $Y_1$ has simple  OSz--homology 
at $\t _1$. Then for every spin$^c$
structure $\s\in Spin^c (Y)$ with $\s \vert _{Y_1}=\t _1$
the maps $F^+_{W,\s}$ and ${\widehat {F}}_{W,\s}$
vanish.
\end{cor}

\begin{proof}
Proposition~\ref{p:cobordism} implies the vanishing of $F^{\infty}_{W,\s}$. 
Combining this with the assumption that $Y_1$ has simple  OSz--homology 
at $\t _1$,
the fact that  
\[
F^+_{W,\s}\circ \varphi_{(Y_1,\s |{Y_1})} = 
\varphi_{(Y_2,\s |_{Y_2})} \circ F^{\infty} _{W,\s}
\]
immediately implies the vanishing of $F^+_{W,\s}$. The vanishing of
${\widehat {F}}_{W,\s}$ now follows from the naturality of the exact
sequence connecting the groups $\hf$ and $HF^+$.
\end{proof}

Examples of OSz--simple manifolds are provided by certain
torus bundles over $S^1$.

\begin{prop}\label{p:ell-par}
A torus bundle $Y\to S^1$ with elliptic or parabolic monodromy $A\in
SL(2, \bfz )$ (that is, $\vert tr (A)\vert <2$ or $\vert tr (A)\vert
=2$) is OSz--simple.
\end{prop}

\begin{proof}
Suppose first that $Y$ has elliptic monodromy. By the classification
of torus bundles over $S^1$ (see e.g.~\cite{Ha}) it follows that, up
to changing its orientation, $Y$ is the boundary of one of the three
plumbings described in Figure~\ref{f:ell}.
\begin{figure}[ht]
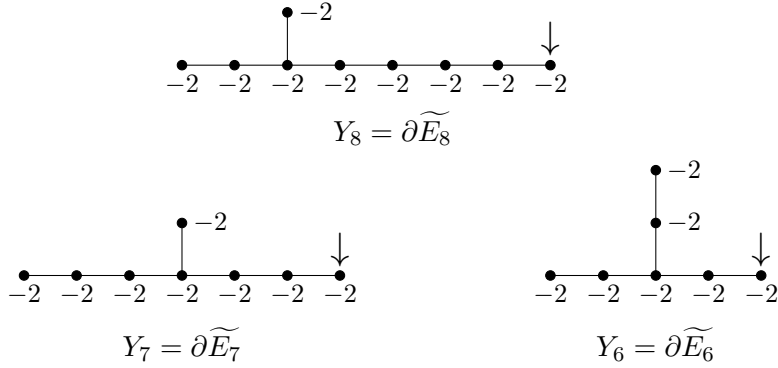

\begin{center}
\setlength{\unitlength}{1mm}
\unitlength=0.7cm
\begin{graph}(18,7.5)(-2,0)
\graphnodesize{0.2}
  \roundnode{m1}(0,2)
  \roundnode{m2}(1,2)
  \roundnode{m3}(2,2)
  \roundnode{m4}(3,2)
  \roundnode{m5}(4,2)
  \roundnode{m6}(5,2)
  \roundnode{m7}(6,2)
  \roundnode{m8}(3,3)

  \edge{m1}{m2}
  \edge{m2}{m3}
  \edge{m3}{m4}
  \edge{m4}{m5}
  \edge{m5}{m6}
  \edge{m6}{m7}
  \edge{m4}{m8}

  \autonodetext{m7}[n]{\Large $\downarrow$}
  \autonodetext{m1}[s]{{\small $-2$}}
  \autonodetext{m2}[s]{{\small $-2$}}
  \autonodetext{m3}[s]{{\small $-2$}}
  \autonodetext{m4}[s]{{\small $-2$}}
  \autonodetext{m5}[s]{{\small $-2$}}
  \autonodetext{m6}[s]{{\small $-2$}}
  \autonodetext{m7}[s]{{\small $-2$}}
  \autonodetext{m8}[e]{{\small $-2$}}

  \freetext(3,0.7){$Y_7=\del\widetilde{E_7}$}

  \roundnode{n1}(10,2)
  \roundnode{n2}(11,2)
  \roundnode{n3}(12,2)
  \roundnode{n4}(13,2)
  \roundnode{n5}(14,2)
  \roundnode{n6}(12,3)
  \roundnode{n7}(12,4)

  \edge{n1}{n2}
  \edge{n2}{n3}
  \edge{n3}{n4}
  \edge{n4}{n5}
  \edge{n3}{n6}
  \edge{n6}{n7}

  \autonodetext{n5}[n]{\Large $\downarrow$}
  \autonodetext{n1}[s]{{\small $-2$}}
  \autonodetext{n2}[s]{{\small $-2$}}
  \autonodetext{n3}[s]{{\small $-2$}}
  \autonodetext{n4}[s]{{\small $-2$}}
  \autonodetext{n5}[s]{{\small $-2$}}
  \autonodetext{n6}[e]{{\small $-2$}}
  \autonodetext{n7}[e]{{\small $-2$}}

  \freetext(12,0.7){$Y_6=\del\widetilde{E_6}$}

  \roundnode{p1}(3,6)
  \roundnode{p2}(4,6)
  \roundnode{p3}(5,6)
  \roundnode{p4}(6,6)
  \roundnode{p5}(7,6)
  \roundnode{p6}(8,6)
  \roundnode{p7}(9,6)
  \roundnode{p8}(10,6)
  \roundnode{p10}(5,7)

  \edge{p1}{p2}
  \edge{p2}{p3}
  \edge{p3}{p4}
  \edge{p4}{p5}
  \edge{p5}{p6}
  \edge{p6}{p7}
  \edge{p7}{p8}
  \edge{p3}{p10}

  \autonodetext{p1}[s]{{\small $-2$}}
  \autonodetext{p2}[s]{{\small $-2$}}
  \autonodetext{p3}[s]{{\small $-2$}}
  \autonodetext{p4}[s]{{\small $-2$}}
  \autonodetext{p5}[s]{{\small $-2$}}
  \autonodetext{p6}[s]{{\small $-2$}}
  \autonodetext{p7}[s]{{\small $-2$}}
  \autonodetext{p8}[s]{{\small $-2$}}
  \autonodetext{p10}[e]{{\small $-2$}}
  \autonodetext{p8}[n]{\Large $\downarrow$}

  \freetext(7,4.8){$Y_8=\del\widetilde{E_8}$}

\end{graph}
\end{center}
\caption{\quad Torus bundles with elliptic monodromy}
\label{f:ell}
\end{figure}
In fact, these plumbings are regular neighbourhoods of the elliptic
singular fibers ${\tilde {E_6}}$, ${\tilde {E_7}}$ and ${\tilde
{E_8}}$, cf.~\cite{HKK}. It is easy to check that by deleting the
vertices indicated by the arrows one gets the 3--manifolds
$S^3_{i-9}(K)$, where $K$ denotes the left--handed trefoil knot
($i=6,7,8$). On the other hand, by assigning weight $(-1)$ to the
vertices indicated by the arrows, we get the lens spaces $L(9-i, 1)$
($i=6,7,8$). Since lens spaces, and all $r$--surgeries on $K$ with
$r\leq -1$ are $L$--spaces~\cite[\S 3]{OSzF2},~\cite[Lemma~7.12 and \S
8]{abs}, the surgery exact triangle \cite{OSzF2} for the $\hf$--theory
implies that $\rank \hf (Y_i)\leq 2(9-i)$ ($i=6,7,8)$. Since
\[
H_1(Y_i; \bfz ) = \bfz \oplus \bfz / (9-i)\bfz \qquad (i=6,7,8),
\] 
and for a torsion spin$^c$ structure $\t \in Spin ^c(Y)$ with
$b_1(Y)=1$ we have $\rank \hf (Y, \t )\geq 2$ (cf. \cite{LSpacific}), we
conclude that $\rank \hf (Y_i)=2(9-i)$, which, in view
of~\cite[Proposition~2.2]{LSpacific} shows that each $Y_i$ ($i=6,7,8$)
is OSz--simple.  (For $Y_8$ the same fact is proved
in~\cite[Section~8.1]{abs}.)

Let us now consider the case of parabolic monodromy.  By the
classification of torus bundles (cf.~e.g.~\cite{Ha}), we know that $Y$
is either a circle bundle over a torus or a Klein bottle, or it is
diffeomorphic to the Seifert fibered 3--manifold described by the
diagram of Figure~\ref{f:seif}.

\begin{figure}[ht]
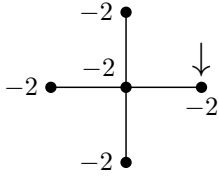

\begin{center}
\setlength{\unitlength}{1mm}
\unitlength=1cm
\begin{graph}(7,2.2)(0,0)

\graphnodesize{0.15}

  \roundnode{n1}(3,1)
  \roundnode{n2}(4,1)
  \roundnode{n3}(5,1)
  \roundnode{n4}(4,0)
  \roundnode{n5}(4,2)

  \edge{n1}{n2}
  \edge{n2}{n3}
  \edge{n2}{n4}
  \edge{n2}{n5}

  \autonodetext{n1}[w]{{\small $-2$}}
  \autonodetext{n2}[nw]{{\small $-2$}}
  \autonodetext{n3}[s]{{\small $-2$}}
  \autonodetext{n4}[w]{{\small $-2$}}
  \autonodetext{n5}[w]{{\small $-2$}}
  \autonodetext{n3}[n]{\Large $\downarrow$}

\end{graph}
\end{center}
\caption{The Seifert fibered manifold
$M(0;\frac{1}{2},\frac{1}{2},-\frac{1}{2},-\frac{1}{2})$}
\label{f:seif}
\end{figure}

The trivial circle bundle $T^3$ over the 2--torus is OSz--simple
by~\cite[Section~8.4]{abs}.  If the circle bundle $Y_n \to T^2$ has
Euler number $n$ then $Y_n$ is diffeomorphic to $M\{ 0, 0, n\}$ of
\cite[Subsection~8.2]{abs}. For $n=1$ the Ozsv\'ath--Szab\'o homology
group $\hf(Y_1)$ is shown in the proof of \cite[Proposition~8.4]{abs}
to be $\Z^4$, verifying the statement.  For $n>1$ we can proceed by a
simple induction on $n$: By the surgery triangle written for the
$n$--framed unknot in the surgery diagram for $M\{ 0, 0, n\}$
described in \cite[Subsection~8.2]{abs} we get
\[
\rank \hf (Y_{n+1} )\leq \rank \hf (S^1\times  S^2\# S^1\times  S^2)+ \rank
\hf (Y_n)=4+4n=4(n+1).
\]
On the other hand, 
\[
H_1(Y_{n+1}; \Z ) \cong \Z^2 \oplus \Z/(n+1)\Z
\]
implies that $\rank \hf (Y_{n+1})\geq 4(n+1)$, hence $Y_n $ is OSz--simple for
$n \geq 0$.  Since by reversing the orientation if necessary, we may assume
the Euler number $n$ to be positive, we conclude the proof for circle bundles
over $T^2$.

Circle bundles over the Klein bottle $K$ can be handled similarly. A surgery
description of such a 3--manifold $Z_n$ with Euler number $n$ is given by
Figure~\ref{f:klein} (cf.~\cite[Figure~6.4 with $k=0$ and $l=2$]{GS}).
\begin{figure}[ht]
\begin{center}
\psfrag{n}{$n+4$}
\psfrag{0}{$0$}
\psfrag{2}{$2$}
\includegraphics[height=2.5cm]{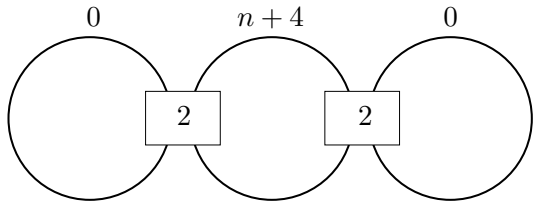}
\end{center}
\caption{\quad Circle bundle over the Klein bottle with Euler number
$n$}
\label{f:klein}
\end{figure}
Simple Kirby calculus shows that this diagram provides the same
3--manifold as the plumbing of Figure~\ref{f:villa}(a). 
\begin{figure}[ht]
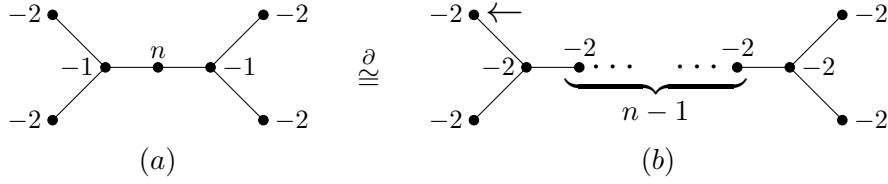

\begin{center}
\setlength{\unitlength}{1mm}
\unitlength=0.7cm
\begin{graph}(15,2.8)(0,-0.5)
\graphnodesize{0.2}

  \roundnode{m1}(0,0)
  \roundnode{m2}(0,2)
  \roundnode{m3}(1,1)
  \roundnode{m4}(2,1)
  \roundnode{m5}(3,1)
  \roundnode{m6}(4,0)
  \roundnode{m7}(4,2)
  \edge{m1}{m3}
  \edge{m2}{m3}
  \edge{m3}{m4}
  \edge{m4}{m5}
  \edge{m5}{m6}
  \edge{m5}{m7}

  \autonodetext{m1}[w]{{\small $-2$}}
  \autonodetext{m2}[w]{{\small $-2$}}
  \autonodetext{m3}[w]{{\small $-1$}}
  \autonodetext{m4}[n]{{\small $n$}}
  \autonodetext{m5}[e]{{\small $-1$}}
  \autonodetext{m6}[e]{{\small $-2$}}
  \autonodetext{m7}[e]{{\small $-2$}}

 \freetext(2,-0.8){$(a)$}

  \freetext(6,1){$\stackrel{\del}{\cong}$}
 
\freetext(11.45,0.37)
{$\underset{{\textstyle n-1}}{\underbrace{\hspace{2.4cm}}}$}

  \roundnode{n1}(8,0)
  \roundnode{n2}(8,2)
  \roundnode{n3}(9,1)
  \roundnode{n4}(10,1)
  \roundnode{n5}(13,1)
  \roundnode{n6}(14,1)
  \roundnode{n7}(15,0)
  \roundnode{n8}(15,2)

  \edge{n1}{n3}
  \edge{n2}{n3}
  \edge{n3}{n4}
  \edge{n5}{n6}
  \edge{n6}{n7}
  \edge{n6}{n8}

  \autonodetext{n1}[w]{{\small $-2$}}
  \autonodetext{n2}[w]{{\small $-2$}}
  \autonodetext{n3}[w]{{\small $-2$}}
  \autonodetext{n4}[n]{{\small $-2$}}
  \autonodetext{n4}[e]{{\Large $\cdots$}}
  \autonodetext{n5}[w]{{\Large $\cdots$}}
  \autonodetext{n5}[n]{{\small $-2$}}
  \autonodetext{n6}[e]{{\small $-2$}}
  \autonodetext{n7}[e]{{\small $-2$}}
  \autonodetext{n8}[e]{{\small $-2$}}
  \autonodetext{n2}[e]{\Large $\leftarrow$}
 \freetext(11.5,-0.8){$(b)$}

\end{graph}
\end{center}
\caption{\quad Alternative plumbing diagrams for 
circle bundle over the Klein bottle}
\label{f:villa}
\end{figure}

For $n>0$ this is equivalent to the plumbing of
Figure~\ref{f:villa}(b), and for $n=0$ (after turning the diagram into
a surgery picture and sliding one $(-1)$--circle over the other and
cancelling the 0--framed unknot against the $(-1)$--circle) we get
that Figure~\ref{f:villa}(a) gives the same 3--manifold as
Figure~\ref{f:seif}.  As before, we can assume that $n\geq 0$ by
possibly reversing orientation.  Consider the surgery exact sequence
for the vertex indicated by the arrow in Figure~\ref{f:villa}(b) (and
Figure~\ref{f:seif} for $n=0$). Notice that the two other manifolds in
the surgery triangle are both $L$--spaces: one is diffeomorphic to the
link $L_{n+4}$ of the $D_{n+4}$ singularity, while the other to $L_n$
for $n\geq 4$, to $L(4,3)$ for $n=3$, to $L(2,1)\# L(2,1)$ for $n=2$,
to $L(4,1)$ for $n=1$ and finally to $-L_4$ for $n=0$. Since $L_n$ is
well--known to have elliptic geometry,
by~\cite[Proposition~2.3]{OSzlens} it is an $L$--space. Thus, since
$|H_1(L_n;\Z)|=4$, we have 
\[
\hf (L_n)\cong \hf (L(4,i))\cong \hf (L(2,1)\# L(2,1))
\cong \Z ^4.
\]
This implies that $\rank \hf (Z_n)\leq 8$.  On the other hand,
$H_1(Z_n ; \Z )$ is either $\Z \oplus \Z / 2\Z \oplus \Z/2\Z$ or $\Z
\oplus \Z / 4\Z$ (depending on the parity of $n$), hence we conclude
that $\hf (Z_n )=\Z ^8$, verifying the statement.
\end{proof}

\begin{rem}
Notice that torus bundles with elliptic monodromies are
boundaries of neighbourhoods of type $II, II^*, III, III^*, IV, IV^*$
fibers in elliptic fibrations \cite{HKK}. Torus bundles with parabolic
monodromies
can be regarded (up to orientation) as boundaries of neighbourhoods of
elliptic 
$I_n$--fibers (when the torus bundle is a circle bundle over $T^2$, $n \geq
1$)
and of elliptic $I_n ^*$--fibers (which are $S^1$--fibrations over the Klein
bottle, $n \geq 0$), cf. \cite{HKK}.
\end{rem}

Further examples of OSz--simple 3--manifolds are provided by certain
Seifert fibered 3--manifolds.  If $Y$ is a Seifert fibered 3--manifold
over $S^2$ with nonnegative background Chern number then $-Y$ is the
boundary of a starshaped plumbing with no bad vertices (in the sense
of \cite{OSzplum}). By \cite{OSzplum} this implies that $-Y$, and
therefore $Y$ is OSz--simple. If $Y$ is a Seifert fibered 3--manifold
over an orientable base and the background Chern number is large
enough then $Y$ is OSz--simple, provided we use $\Z/2\Z$--coefficients
in the definition of the Ozsv\'ath--Szab\'o homology groups.  Most
probably the statement holds true for $\Z$--coefficients as well, but
since the computational tool we will use in the proof has been
verified in \cite{OSzrat} with $\Z/2\Z$--coefficients, we restrict our
attention to this case.  Before stating the result we need to fix our
notations on Seifert fibered spaces. We do this following the
conventions of~\cite{OSzrat}. We say that a Seifert fibered
3--manifold $Y$ over a genus $g$ surface has~\emph{Seifert invariants}
$(a,\frac{r_1}{q_1},\ldots,\frac{r_n}{q_n})$ if the Seifert fibration
on $Y$ is obtained in the canonical way by performing $(-\frac{q_1
}{r_1})$--, $\ldots, (-\frac{q_n }{r_n})$--surgeries along $n$ fibers
of the circle bundle $Y_{g,a}\to \Sigma_g$ over an orientable
genus--$g$ surface $\Sigma_g$ with Euler number $e(Y_{g,a})=a$.

\begin{prop}\label{p:seif}
Let $Y$ be a Seifert fibered 3--manifold over a genus $g$ surface 
with Seifert invariants $(a,\frac{r_1}{q_1},\ldots,\frac{r_n}{q_n})$.  
If $a> 2g$ and we consider Ozsv\'ath--Szab\'o homology groups with 
$\Z/ 2\Z$--coefficients, then $Y$ is OSz--simple.
\end{prop}

\begin{proof} According to \cite[Theorem 10.1]{OSzrat} 
we only need to check that the function $h_t\co\Z\to\Z$ has a unique
local minimum for every $t$ satisfying $-g\leq t \leq g$. By
definition,
\[
h_t(s) = 
\begin{cases}
\sum_{i=0}^{s-1} \de_t(i)\qquad\text{if $s> 0$}\\
0 \qquad\qquad\quad\quad\text{if $s=0$}\\
-\sum_{i=s}^{-1} \de_t(i)\quad\text{if $s<0$},
\end{cases}
\] 
where 
\begin{equation}\label{e:deltat}
\delta_t (s)=(-1)^{s+1}t + \xi_0 +a\cdot s +
\sum_{i=1}^n \lfloor\frac{\xi_i +r_i \cdot s}{q_i}\rfloor.
\end{equation}
Therefore it suffices to show that, for every $-g\leq t\leq g$, the
function $\de_t$ changes sign only once, that is,
\[
\delta _t(s)>0\quad\text{implies}\quad\delta _t (s+1)>0.
\]
Notice that we have a freedom in choosing $\xi _0$ (and \cite[Theorem
  10.1]{OSzrat} shows, in particular, that different choices giving
the same spin$^c$ structure yield the same Ozsv\'ath--Szab\'o homology
groups). Since we are only concerned with torsion spin$^c$ structures,
we can fix $\xi _0$ to be arbitrarily large in absolute value and
negative. It then follows from Formula~\eqref{e:deltat} that we can
assume $\delta _t(s)<0$ for $s\leq 0$. Now suppose that $-g\leq t\leq
g$, $s>0$ and $\delta _t (s)>0$. To finish the proof it clearly
suffices to verify that $\de_t(s)<\de_t(s+1)$. Since
\[
\sum _{i=1}^n \lfloor\frac{\xi_i +r_i \cdot s}{q_i}\rfloor < 
\sum_{i=1}^n \lfloor\frac{\xi _i + r_i \cdot (s+1)}{q_i}\rfloor
\]
and $2\vert t \vert \leq 2g <a$, we have 
\begin{multline*}
\delta_t (s)=(-1)^{s+1}t + \xi_0 + a\cdot s + 
\sum_{i=1}^n \lfloor\frac{\xi_i +r_i \cdot s}{q_i}\rfloor <\\
- (-1)^{s+1}t + \xi_0 + a\cdot s + a + 
\sum _{i=1}^n \lfloor\frac{\xi_i + r_i \cdot (s+1)}{q_i}\rfloor
= \delta _t (s+1).
\end{multline*}
\end{proof}

\sh{Contact Ozsv\'ath--Szab\'o invariants}

A contact structure $\xi $ on $Y$ determines an element $c(Y,\xi )\in
\hf (-Y, \t _{\xi })$ (and similarly in $HF ^+(-Y , \t _{\xi })$)
up to sign, which has the following crucial
properties~\cite{OSzcont, OSzgen}:
\begin{itemize}
\item $\pm c(Y, \xi )$ is an isotopy invariant of the contact 3--manifold
$(Y, \xi )$;
\item $c(Y, \xi )=0$ if the contact structure $\xi $ is overtwisted;
\item $c(Y, \xi )\neq 0$ if $(Y, \xi )$ is strongly  fillable;
\item if $(Y_L, \xi _L)$ is given by contact $(-1)$--surgery along 
the Legendrian knot $L\subset (Y, \xi)$, inducing the Stein cobordism $X$
with spin$^c$ structure $\s _0$ then by ${\overline {X}}$ 
denoting $X$ when turned upside down we have
\[
F_{{\overline {X}}, \s _0}(c(Y_L,\xi _L))=c(Y, \xi )\qquad {\mbox
{and}} \qquad F_{{\overline {X}}, \s}(c(Y_L,\xi _L))=0
\]
for all other spin$^c$ structures $\s \neq \s _0$, \cite{Gh1, Plam}.
\end{itemize}

\begin{rem}
The above statements hold true using both $\Z$ and $\Z/2\Z$
coefficients. Ozsv\'ath and Szab\'o defined twisted versions of their
contact invariants in such a way that every weakly fillable contact
structure admits a nontrivial twisted contact Ozsv\'ath--Szab\'o
invariant for some appropriate twisting. In this paper, however, we
concentrate on \emph{untwisted} invariants.
\end{rem}

\section{Auxiliary results}
\label{s:third}

In this section we establish two auxiliary results which will be 
used in the proofs of the next section.

Suppose that $T^2\times [0,1]$ is embedded into a 3--manifold
$Y$. Consider the tori $T_i=T^2\times \{ t_i \}$ with $0<t_1 < t_2 <
\ldots < t_k<0$ and for every $i=1,\ldots , k$ let
$\{C^j_i\}_{j=1}^{s_i}\subset T_i$ be a finite collection of parallel
and disjoint simple closed curves. Perform 3--dimensional Dehn surgery
along each $C^j_i$ with framing $-1$ with respect to the framing
induced by the torus $T_i$, and call $Y'$ the resulting
3--manifold. In the following $D_C$ will denote the right--handed Dehn
twist along the curve $C\subset T^2$ in the mapping class group $\Gamma
_1$ of the torus $T^2$.

\begin{prop}\label{p:order}
The 3--manifold $Y'$ is obtained from $Y$ by cutting it along
$T^2\times \{ 0 \}$ and regluing via the diffeomorphism
\[
D^{s_k}_{C_k}\circ D^{s_{k-1}}_{C_{k-1}}\circ \ldots \circ D^{s_1}_{C_1}.
\]
\end{prop}

\begin{proof}
It is an easy exercise to check that the surgery along each $C^j_i$
results in cutting $Y$ along $T\x\{t_i\}$ and regluing with the map
$D_{C^j_i}$. To prove the statement we only need to check that
performing all the surgeries is equivalent to cutting and regluing via
the composition of diffeomorphisms in the order stated. In order to
see this, modulo an easy induction argument it suffices to show that
if $F$, $G$ and $H$ are closed, oriented surfaces and $\varphi\co F\to
G$, $\psi\co G\to H$ are orientation--preserving diffeomorphisms, then the
two quotients
\[
\frac
{(F\x [0,\frac13])\sqcup (G\x [\frac13,\frac23])\sqcup 
(H\x [\frac23,1])} 
{(x,\frac13)\sim (\varphi(x),\frac13),\ 
(y,\frac23)\sim (\psi(y),\frac23)}
\]
and 
\[
\frac
{(F\x [0,\frac23])\sqcup (H\x [\frac23,1])}
{(x,\frac23)\sim ((\psi\circ\varphi)(x),\frac23)}
\]
are orientation--preserving diffeomorphic. In fact, an
orientation--preserving diffeomorphism is induced by the map
\[
(\id_F\x \id_{[0,\frac13]})\sqcup(\varphi^{-1}\x\id_{[\frac13,\frac23]})
\sqcup(\id_H\x\id_{[\frac23,1]}).
\]
\end{proof}

Let us now fix an identification of $T^2 \times \{ 0\}$ with $\R ^2
/\Z ^2$, and let $a$ and $b$ denote the linear curves with slopes $0$
and $\infty$, respectively, obtained by mapping the coordinate axes of
$\R^2$ to $\R ^2 /\Z ^2$. For short, let us also denote by $a$ and $b$
the right--handed Dehn twists along the curves $a$ and $b$. It is a
well--known fact that $a$ and $b$ generate the mapping class group
$\Gamma _1$, which has presentation
\[
\Gamma _1 = \langle a, b \mid aba=bab, \ (ab)^6=1 \rangle.
\]
Using the relation $aba=bab$ it easily follows from $(ab)^6=1$ that
$(a^3b)^3=(b^3a)^3=1$. Consider the element
\begin{equation}\label{e:element}
\ga = (a^3b)^3 b = a^3 b a^3 b a^3 b^2
\end{equation}
in $\Gamma_1$ viewed as a product of six factors, each of which is a
power of either $a$ or $b$. Let $0<t_1<\cdots < t_6<1$, and consider
simple closed curves $C_i\subset T^2\x\{t_i\}$ with $C_i$ isotopic to
$b$ for $i$ odd and to $a$ for $i$ even. By adding the right number of
parallel copies of the same curve on each torus $T^2\x\{t_i\}$ we can
ensure that the diffeomorphism associated via
Proposition~\ref{p:order} to performing $(-1)$--surgery along each of
the curves is the above elemenent $\ga\in\Ga_1$. Attach 4--dimensional
2--handles along the above knots with framing $(-1)$ with respect to
the surface framings induced by the tori $T^2\x\{t_i\}$, and denote
the resulting 4--dimensional cobordism built on $Y$ by $W$.

\begin{prop}\label{p:sign}
The 4--dimensional cobordism $W$ defined above satisfies $b_2^+(W)>0$.
\end{prop}

\begin{proof}
Consider the first six $a$--curves $C_i$ and their corresponding
2--handles. By sliding each of the first five 2--handles over the next
one it is easy to see that $W$ contains the 4--manifold obtained by
attaching 2--handles to a chain of five $(-2)$--framed unknots contained
in a 3--ball inside $Y$. Now consider the four $b$--curves and slide
their corresponding 2--handles in the same way. This gives a framed
link still contained in a 3--ball inside $Y$ with intersection graph
given by Figure~\ref{f:plum}.
\begin{figure}[ht]
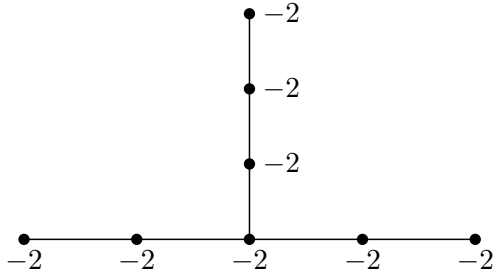

\begin{center}
\setlength{\unitlength}{1mm}
\unitlength=1cm
\begin{graph}(8,3)(0,0)
\graphnodesize{0.15}
  \roundnode{m1}(1,0)
  \roundnode{m2}(2.5,0)
  \roundnode{m3}(4,0)
  \roundnode{m4}(5.5,0)
  \roundnode{m5}(7,0)
  \roundnode{top1}(4,1)
  \roundnode{top2}(4,2)
  \roundnode{top3}(4,3)

  \edge{m1}{m2}
  \edge{m2}{m3}
  \edge{m3}{m4}
  \edge{m4}{m5}
  \edge{m3}{top1}
  \edge{top1}{top2}
  \edge{top2}{top3}

  \autonodetext{m1}[s]{$-2$}
  \autonodetext{m2}[s]{$-2$}
  \autonodetext{m3}[s]{$-2$}
  \autonodetext{m4}[s]{$-2$}
  \autonodetext{m5}[s]{$-2$}
  \autonodetext{top1}[e]{$-2$}
  \autonodetext{top2}[e]{$-2$}
  \autonodetext{top3}[e]{$-2$}
\end{graph}
\end{center}
\caption{The 4--dimensional plumbing $P\subseteq W$}
\label{f:plum}
\end{figure}

Viewing Figure~\ref{f:plum} as representing a smooth 4--dimensional
plumbing $P\subseteq W$, it is easy to check that the associated $8\x
8$ intersection matrix has determinant $-3$, hence the matrix has at
least one positive eigenvalue. This implies that $b_2^+(P)>0$ and
therefore $b_2^+(W)>0$.
\end{proof}

\section{Proofs of Theorems~\ref{t:main}, \ref{t:bundle}, 
\ref{t:sfs} and~\ref{t:ghiggini+}}
\label{s:fourth}

\sh{Proof of Theorem~\ref{t:main}}

The assumption that $(Y,\xi)$ has torsion at least $2$ implies that
there is a contact embedding $\T_2\hra (Y,\xi)$. Since the germ of a
contact structure around a surface is determined by the induced
characteristic foliation, for some small $\ep>0$ there is a contact
embedding
\begin{equation}\label{e:embedding}
(T^2\times [-\ep,1+\ep], 
\ker(\cos (4\pi z)dx - \sin (4 \pi z)dy)) \hra (Y,\xi).
\end{equation}
Fix an identification of the torus $T^2$ with $\R^2/\Z^2$ such that the
characteristic foliation induced on $T^2\x\{0\}$ has slope $\infty$ and let
$a$ and $b$ be simple closed curves on $T^2$ with slopes $0$ and $\infty$,
respectively.  As before, by abuse of notation, we shall denote by $a$ and $b$ the
elements of the mapping class group of the torus $\Gamma_1$ determined by
positive Dehn twists around the curves $a$ and $b$.  Since in the group
$\Gamma_1$ we have $(a^3b)^3=(b^3a)^3=1$, it follows that
\[
1 = (a^3b)^3(b^3a)^3 = a^3ba^3ba^3b^4ab^3ab^3a,
\]
therefore
\begin{equation}\label{e:word}
1 = a(a^3b)^3(b^3a)^3a^{-1} = a^4ba^3ba^3b^4ab^3ab^3.
\end{equation}
Then conjugating the last word in~\eqref{e:word} by $b^n$ we get the word  
\begin{equation}\label{e:newword}
(b^nab^{-n})^4 b (b^nab^{-n})^3b (b^nab^{-n})^3b^4 (b^nab^{-n})
b^3(b^nab^{-n}) b^3,
\end{equation}
which is easily checked to be a composition of powers of Dehn twists along
simple closed curves of slopes
\begin{equation}\label{e:slopes}
-n, \infty , -n, \infty, -n, \infty , -n, \infty , -n, \infty.
\end{equation}
If $n$ is sufficiently large, by~\eqref{e:embedding} we can locate
inside $(Y,\xi)$ embedded tori whose characteristic foliations are
made of simple closed curves having slopes given
by~\eqref{e:slopes}. In view of Proposition~\ref{p:order} we can
perform Legendrian surgery on a suitable number of parallel Legendrian
curves on such tori, so that the resulting smooth 3--manifold is
obtained by cutting $Y$ along $T$ and regluing via a diffeomorphism
whose isotopy class is specified by the word~\eqref{e:newword}. But by
construction this word represents $1\in\Ga_1$, therefore the resulting
3--manifold is $Y$ again, so the construction yields a Stein cobordism
$W$ from $(Y,\xi)$ to $(Y,\xi')$ for some contact structure
$\xi'$. Moreover, we know that $c(Y,\xi)={\widehat {F}}_{{\overline
{W}}, \s _0}(c(Y, \xi'))$. Since the word~\eqref{e:word} contains the
word given by~\eqref{e:element}, by Proposition~\ref{p:sign} we have
$b_2^+(W)>0$. Therefore if $Y$ is OSz--simple then
Corollary~\ref{c:vanish} implies that the map ${\widehat
{F}}_{{\overline {W}}, \s _0}$ vanishes, verifying that $c(Y,\xi)=0$.

We shall now prove the second part of the statement assuming
Theorem~\ref{t:reducetorsion}, which will be proved in
Section~\ref{s:appendix}.  Suppose that $b_1(Y)\leq 1$ and
$\Tor(Y,\xi)>0$, that is, we have a contact embedding $\T_1\hra
(Y,\xi)$. In Section~\ref{s:appendix} we will show that one can build
a Stein cobordism $W_1$ from $(Y,\xi)$ to a contact 3--manifold
$(Y,\xi')$ using the word $(b^3a)^3$, in such a way that $(Y,\xi)$ is
obtained from $(Y,\xi')$ by a Lutz modification. It follows
by~\cite{Colin1} that the contact structures $\xi $ and $\xi'$ are
homotopic as 2--plane fields, and therefore if $c(Y, \xi )\neq 0$ the
two invariants $c(Y,\xi)$ and $c(Y,\xi')$ are elements of the same
Ozsv\'ath--Szab\'o group $\hf_d(-Y,\t_{\xi })$.  Moreover, since $Y$
is OSz--simple and $b_1(Y)\leq 1$, this group is of rank $1$
(cf.~\cite[Definition~4.9]{abs} when $H_1(Y; \bfz )\cong \bfz$, the
remark on \cite[page~250]{abs} in general, and also
\cite[Proposition~2.2]{LSpacific}). Thus, the restriction of the map
$F_{{\overline W_1}, \s _0}$ to $\hf_d(-Y, \t_{\xi })$ is
multiplication by some $k\in\Z$. Now observe that the smooth cobordism
$W_2=W_1\circ W_1$ obtained by ``composing'' $W_1$ with itself can be
constructed using the word $(b^3 a)^6$, which contains (up to
conjugation) the word given by~\eqref{e:element}. Therefore by
Proposition~\ref{p:sign} $b_2^+(W_2)>0$. As before, by
Corollary~\ref{c:vanish} this implies $F^{\infty }_{{\overline {W_2}},
\s }=0$ for any $\s \in Spin ^c ({\overline {W_2}})$, and so by
Theorem~\ref{t:composition}
\[
(F_{{\overline {W_1}}, \s
_0}\circ F_{{\overline {W_1}}, \s _0}) (c(Y, \xi '))=
\sum \pm F_{{\overline {W_2}}, \s}(c(Y, \xi '))=0,
\]
and it follows that $k=0$. This shows that $c(Y, \xi )=F_{{\overline
{W_1}}, \s _0} (c(Y, \xi '))= 0$, concluding the proof.

\sh{Proof of Theorem~\ref{t:bundle}} 

If the monodromy is trivial then $Y$ is the 3--torus $T^3$. Suppose
that $\xi$ is a contact structure on $T^3$ with $\Tor (T^3, \xi)>0$.
By the classification of tight contact structures on
$T^3$~\cite{Ka}, up to applying a diffeomorphism of $T^3$ we may
assume that there is a contact embedding $\T_1\hra (T^3, \xi)$ such
that $T^2\x\{0\}\subseteq\T_1$ maps to $T^2\x\{s\}\subseteq T^2\x
S^1=T^3$ for some $s\in S^1$. Fix an identification of $T^2\x\{s\}$
with $\R^2/\Z^2$, and denote, as before, by $a$ and $b$, respectively,
the right--handed Dehn twists along simple closed curves with slopes
$0$ and $\infty$. Arguing as in the proof of Theorem~\ref{t:main} we
can use the word $(a^3b)^3b = b$ in the mapping class group to build a
Stein cobordism $W$ from $(T^3,\xi)$ to $(Y_1,\xi')$, where $Y_1$ is a
torus bundle over $S^1$ with monodromy $b$. By
Proposition~\ref{p:ell-par} the 3--manifold $Y_1$ is OSz--simple and
by Proposition~\ref{p:sign} we have $b_2^+(W)>0$, therefore by
Corollary~\ref{c:vanish} it follows that $c(T^3, \xi) = F_{{\overline
W}, \s _0}(c(Y_1,\xi')) = 0$.

The proof of the statement when $|\tr(A)|<2$ or $|\tr(A)|=2$ follows
from Proposition~\ref{p:ell-par} combined with Theorem~\ref{t:main}.

Now suppose that $|\tr(A)|> 2$. In this case any incompressible torus
is isotopic to the fiber of the fibration~\cite[Lemma~2.7]{Ha}. Let
the monodromy of the fibration be denoted by $A\in SL_2 (\Z ) \cong
\Gamma _1$ and fix a decomposition of $A^{-1}$ into the product of
right--handed Dehn twists. Since $\T_1$ contains tori with linear
characteristic foliations of any rational slope, if $\T_n\hra (Y,\xi)$
with $n$ sufficiently large, by performing suitable Legendrian
surgeries as before we can construct a Stein cobordism $W_A$ from
$(Y,\xi)$ to $(T^3,\xi_A)$ for some contact structure $\xi_A$.
Moreover, up to choosing a larger $n$ we may assume that $\Tor (T^3,
\xi_A)>0$, and therefore by the first part of the proof $c(T^3,
\xi_A)=0$. It follows that $c(Y, \xi)=F_{{\overline {W_A}}, \s
_0}(c(T^3, \xi _A))= 0$.  Notice that in this way a bound for the
optimal $n_Y$ can be easily deduced from the decomposition of $A^{-1}$
into the product of right--handed Dehn twists. (This bound is still
far from the value $n_Y=0$ predicted by Conjecture~\ref{c:torsion}.)

\sh{Proof of Theorem~\ref{t:sfs}}

By e.g.~\cite[page 30]{Ha}, unless $Y$ is an elliptic or parabolic
torus bundle over $S^1$, an incompressible torus $T\hookrightarrow Y$
can be isotoped to be the union of regular fibers of the Seifert
fibration. Therefore, in view of Theorem~\ref{t:bundle} we may assume
that $T$ consists of regular fibers. By assumption there is a contact
embedding $\T_n\hra (Y,\xi)$ with $n>2$, and we can write
$\T_n=\T_{n-1}\cup \T_1$. Since $\T_1$ contains tori with linear
characteristic foliations with any rational slope, we may assume that
for every integer $k\geq 0$ one of those tori contains $k$ disjoint
Legendrian knots $L_1,\ldots, L_k$ each of which is smoothly isotopic
to a regular fiber of the fibration, and such that the contact
framings and the framings induced by the torus (i.e. by the fibration)
coincide.  Performing Legendrian surgeries along $L_1, \ldots, L_k$
gives a Stein cobordism $W_L$ from $(Y,\xi)$ to a contact Seifert
fibered 3--manifold $(Y',\xi')$ such that, when choosing $k$
sufficiently large, $Y'$ is a Seifert fibered 3--manifold over an
orientable base with background Chern number sufficiently high. By
Proposition~\ref{p:seif} the 3--manifold $Y'$ is OSz--simple
(with $\Z/2\Z$--coefficients). By the
construction we have $\T_{n-1}\hra (Y',\xi')$ and by assumption
$n-1>1$, therefore Theorem~\ref{t:main} implies that $c(Y', \xi
')=0$. Thus $c(Y, \xi )= F_{{\overline {W_L}}, \s _0}(c(Y', \xi ')) =
0$.

\sh{Proof of Theorem~\ref{t:ghiggini+}}

 In this proof we
assume familiarity with the work of Ghiggini~\cite{Gh1, Gh2}. Ghiggini
considers a family $\{\ze_i\}$ (denoted $\{\eta_i\}$ in~\cite{Gh1,
Gh2}) of contact structures on $-\Si$, where the index $i$
varies in the set
\[
{\cal P}_n = \{-n+1,-n+3,\ldots,n-3,n-1\}.
\]
Let $(M_0,\xi_1)$ denote the Stein fillable contact 3--manifold
obtained by Legendrian surgery on the Legendrian right--handed trefoil
with $\tb=+1$ in $(S^3,\xi_{st})$. Each contact structure $\ze_i$ is
constructed by performing Legendrian surgery along a Legendrian knot
inside $(M_0,\xi_1)$. Ghiggini also considers a different tight
contact structure $\xi_n$ on $M_0$, and defines a contact structure
$\eta_0$ on $-\Si$ by Legendrian surgery along a Legendrian knot in
$(M_0,\xi_n)$. Denoting by $\overline\xi$ the contact structure $\xi$
with reversed orientation, Ghiggini shows that $\overline\ze_i$ is
isotopic to $\ze_{-i}$ for every $i\in{\cal P}_n$ and $\eta_0$ is
isotopic to $\overline\eta_0$. All of the above holds regardless of
the parity of $n$. Since the statement has been already proved
in~\cite{Gh2} for every even $n$, from now on we shall assume $n$ odd.
Arguing as in~\cite[Lemma~4.4]{Gh2} and~\cite[Proof of
Theorem~2.4]{Gh2}, it follows that
\begin{equation}\label{e:combinations}
c^+(\eta_0) =
\al_0 c^+(\ze_0) + \sum_{i\in{\cal P}\setminus\{0\}} \al_i (c^+(\ze_i)
+ c^+(\ze_{-i})),
\end{equation}
where we may assume $\al_i\in\{0,1\}$ (it suffices to work with
$\Z/2\Z$--coefficients).  Recall that each of the Legendrian surgeries
from $(M_0,\xi_1)$ to $(-\Si,\ze_i)$ as well as the Legendrian surgery
from $(M_0,\xi_n)$ to $(-\Si,\eta_0)$ induce the same smooth
4--dimensional cobordism $V$. Thus, arguing as in~\cite[Proof of
Theorem~2.4]{Gh2} we have
\[
F^+_{\overline V} (c^+(\ze_i) + c^+(\ze_{-i})) = 0 \quad (\bmod \ 2)
\]
for every $i\in{\cal P}\setminus\{0\}$. Therefore, in view of
Equation~\eqref{e:combinations} we have
\[
\al_0 c^+(\xi_1) = \al_0 F^+_{\overline V} (c^+(\ze_0)) = 
F^+_{\overline V} (c^+(\eta_0)) = c(\xi_n) = 0,
\]
where the last equality follows from Theorem~\ref{t:bundle} because
$M_0$ is a torus bundle with elliptic monodromy and by construction
$\Tor(M_0,\xi_n)\geq n-1>1$ because $n\geq 3$. Since $(M_0,\xi_1)$ is
Stein fillable, we have $c^+(\xi_1)\neq 0$, therefore we conclude that
$\al_0=0$, and from this point on the argument proceeds as
in~\cite[Proof of Theorem~2.4]{Gh2}.

\section{Proof of Theorem~\ref{t:reducetorsion}}\label{s:appendix}

The proof of Theorem~\ref{t:main} relied on the construction of a
particular cobordism $W$ from $Y$ to $Y$ which, provided the contact
structure $\xi $ on $Y$ had torsion $\Tor (Y, \xi )>1$, also supported
a Stein structure.  The chosen Stein cobordism might seem to be
\emph{ad hoc}, but as we explain below, the contact surgery on
$(Y, \xi )$ corresponding to this Stein cobordism has a clear contact
topological interpretation: it is the inverse of a Lutz modification.

In this section we shall assume familiarity with results, notation and
terminology from \cite{Gic, Ho1}. Let $T^2$ be a 2--torus with an
identification $T^2\cong\R^2/\Z^2$. Let $B_0\cong T^2\x [0,1]$ and
$B_1\cong T^2\x [0,1]$ be basic slices of the same sign with boundary
slopes respectively $(s_0,s)$ and $(s,s_1)$. Let $B$ denote $B_0\cup
B_1$, the contact 3--manifold obtained by gluing together (via the
identity map) $B_0$ and $B_1$ along their boundary components of slope
$s$. Let $T\subset B$ be a minimal convex torus parallel to the
boundary having slope $s$, and let $C\subset T$ be one of its
Legendrian divides.

\begin{lem}\label{l:basicslices}
The contact 3--manifold $B'$ obtained from $B$ by Legendrian surgery
along $C$ is isomorphic to the contact 3--manifold $B_0\cup B'_1$,
where $B'_1$ is a basic slice with the same sign as $B_0$ and boundary
slopes $(s,D_C^{-1}(s_1))$.
\end{lem}

\begin{proof}
It is easy to see that since $B$ is the union of two basic slices, it
is contactomorphic to a toric layer sitting inside a neighborhood of a
Legendrian knot in the standard contact 3--sphere $(S^3,\xi_{st})$. It
follows that $B'$ contact embeds into a closed contact 3--manifold
$(Y,\ze)$ given by a Legendrian surgery on $(S^3,\xi_{st})$ in
such a way that the image of any torus in $B'$ parallel to the
boundary bounds a solid torus in $Y$. Since $(Y,\ze)$ is Stein
fillable and hence tight, $B'$ must be both tight and minimally
twisting, otherwise one could easily find an overtwisted disk inside
$(Y,\ze)$. We can choose the identification $T^2\cong\R^2/\Z^2$ so
that $s_0=1$ and $s=0$. Then $s_1=-\frac{1}{n}$ for some integer $n\geq 1$
and the action of $D^{-1}_C$ on $s_1$ is determined by
\[
\begin{pmatrix}
1 & -1 \\
0 & \phantom{-}1
\end{pmatrix}
\begin{pmatrix}
n\\-1
\end{pmatrix}
=
\begin{pmatrix}
n+1\\-1
\end{pmatrix},
\]
which shows that the boundary slopes of $B'$ are $1$ and
$-\frac{1}{n+1}$, therefore $B'$ decomposes as $B_0\cup B'_1$. It
remains to check that the signs of $B'_1$ and $B_0$ are equal. Observe
that $B_1\subset\overline B_1$, where $\overline B_1$ is a basic slice
with boundary slopes $(0,1)$ and $B_0\cup\overline B_1$ embeds in
$(T^3,\ker(\cos(4\pi z)dx - \sin(4\pi z)dy))$, which is symplectically
fillable~\cite{Gt}. By doing Legendrian surgery on $B_0\cup\overline
B_1$ along $C$ and computing as before we see that $B'\subset
B''=B_0\cup\overline B'_1$, where $\overline B'_1$ is a basic slice
with boundary slopes $(0,\infty)$. On the other hand, $B''$ is tight,
minimally twisting and has boundary slopes $(1,\infty)$, hence it is a
basic slice as well. Therefore by Honda's Gluing Theorem~\cite[\S
4.7.4]{Ho1} the sign of $\overline B'_1$ must be the same as the sign
of $B_0$. But the inclusion $B'_1\subset \overline B'_1$ implies,
again by the Gluing Theorem, that the sign of $B'_1$ must be the same
as the sign of $\overline B'_1$. This concludes the proof.
\end{proof}

Suppose that $a,b\in\R$, $a<b$, and define 
\[
\T_n[a,b]=(T^2\times [a,b],
\ker (\cos (2\pi nz)dx - \sin (2\pi nz)dy)).
\]
In this notation, we have $\T_n[0,1]=\T_n$, where $\T_n$ is defined in
Section~\ref{s:intro}. Suppose that $a<c<b$, the characteristic
foliation $\FF$ on the torus $T^2\x\{c\}\subset\T_n[a,b]$ is a union
of simple closed curves, and let $C\subset T^2\x\{c\}$ be such a
closed curve. Then, there is a diffeomorphism
\[
D_C\co T^2\x\{c\}\to T^2\x\{c\}
\]
representing the right--handed Dehn twist along $C$ and such that
$D_C(\FF)=\FF$.

\begin{lem}\label{l:surg=gluing}
The contact 3--manifold obtained from $\T_n[a,b]$ by Legendrian
surgery along $C$ is isomorphic to the contact 3--manifold
\[
\T_n[a,c]\cup_{D_C} \T_n[c,b]
\]
obtained by gluing $\T_n[a,c]$ to $\T_n[c,b]$ via the diffeomorphism
$D_C$.
\end{lem}

\begin{proof}
Suppose that the torus $T=T^2\x\{c\}$ has slope $s$. Then $T$ can be
slightly perturbed to become a convex torus with minimal dividing set
of slope $s$ in such a way that a closed leaf $C$ of the
characteristic foliation becomes a Legendrian
divide~\cite[Lemma~3.4]{Gh0}.  We can choose $c_0\in (a,c)$ and $c_1\in
(c,b)$ so that the tori $T_0=T^2\x\{c_0\}$ and $T_1=T^2\x\{c_1\}$ can
be perturbed to minimal convex tori with boundary slopes $s_0$ and
$s_1$, respectively, making sure at the same time that the resulting
layers $B_0$ between $T_0$ and $T$ and $B_1$ between $T$ and $T_1$ are
both basic slices. Since
\[
\T_n[a,b]\subset (T^2\times\R,
\ker (\cos (2\pi nz)dx - \sin (2\pi nz)dy)),
\]
using the Gluing Theorem as in the proof of Lemma~\ref{l:basicslices}
one can easily check that $B_0$ and $B_1$ must have the same sign.
We have the decomposition 
\[
\T_n[a,b] = N_0\cup B_0\cup B_1\cup N_1,
\]
where each of $N_0$ and $N_1$ is a toric layer with only one convex 
boundary component. In view of Lemma~\ref{l:basicslices}, the result 
of Legendrian surgery along $C$ can be decomposed as 
\[
N_0\cup B_0\cup B'_1\cup_{D_C} N_1,
\]
where $B_0$ and $B'_1$ have the same sign and $B'_1$ is glued to $N_1$
via the diffeomorphism $D_C$. But it is easy to check that this is
exactly the decomposition of the contact 3--manifold
$\T_n[a,c]\cup_{D_C} \T_n[c,b]$ obtained by perturbing $T^2\x\{c_0\}$,
$T^2\x\{c\}$ and $T^2\x\{c_1\}$ to become minimal convex tori.
\end{proof}

\begin{proof}[Proof of Theorem~\ref{t:reducetorsion}]
Arguing as in the proof of Theorem~\ref{t:main} we see that for some
small $\ep>0$ there is a contact embedding
\begin{equation*}
(T^2\times [-\ep,1+\ep], 
\ker(\cos (2\pi nz)dx - \sin (2\pi nz)dy)) \hra (Y,\xi).
\end{equation*}
We can choose $\de>0$, $-\ep<-\de<0$, and an identification
$T^2\cong\R^2/\Z^2$ so that the characteristic foliations on
$T_{-\de}=T^2\x\{-\de\}$ and $T_0=T^2\x\{0\}$ are made of simple
closed curves and have slope, respectively, $0$ and $\infty$. Then up
to reparametrizing the interval $[0,1]$ we may assume that the
characteristic foliations on $T_{1/4}=T^2\x\{1/4\}$,
$T_{1/2}=T^2\x\{1/2\}$ and $T_{3/4}=T^2\x\{3/4\}$ are made of simple
closed curves and have slopes, respectively, $0$, $\infty$ and
$0$. Let
\[
C_{-\de}\subset T_{-\de},\quad
C_{1/4}\subset T_{1/4},\quad
C_{3/4}\subset T_{3/4}
\]
and
\[
D^1_0, D^2_0, D^3_0\subset T_0,\quad 
D^1_{1/2}, D^2_{1/2}, D^3_{1/2}\subset T_{1/2},\quad
D^1_1, D^2_1, D^3_1\subset T_1 
\]
be disjoint closed leaves of the respective characteristic
foliations. Observe that if we perform Legendrian surgery on $(Y,\xi)$
along (the images of) each of the curves $C$'s and $D$'s we obtain a
contact 3--manifold of the form $(Y,\xi')$. In fact, each $C$--curve
has slope $0$, while each $D$--curve has slope $\infty$. Therefore, if
we denote by $A$, respectively $B$, the corresponding Dehn twists up
to isotopy, since in the mapping class group of the torus
\begin{equation}\label{e:relation}
B^3 A B^3 A B^3 A = (B^3 A)^3 = 1,
\end{equation}
if follows from Proposition~\ref{p:order} that the 3--manifold
underlying the result of the Legendrian surgeries is still $Y$. To see
that $(Y,\xi)$ is obtained from $(Y,\xi')$ by a Lutz modification,
observe that the tori $T_{-\de}$, $T_0$, $T_{1/4}$, $T_{1/2}$,
$T_{3/4}$ and $T_1$ induce, for some $\de'>0$, $\ep>\de'>0$, a
decomposition
\[
(T^2\times [-\de,1+\de'], \xi) = 
N_1\cup N_2\cup N_3\cup N_4 \cup N_5\cup N_6,
\]
where the boundary components of $N_6$ have characteristic
foliations of slopes $(\infty,n)$ for some $n\geq 1$. When we perform
the Legendrian surgeries along $D^1_1$, $D^2_1$, $D^3_1$ and
$C_{3/4}$, according to Lemma~\ref{l:surg=gluing} the above
decomposition becomes
\[
N_1\cup N_2\cup N_3\cup N_4\cup_A N_5\cup_{B^3} N_6.
\]
Since 
\[
A^{-1}
\begin{pmatrix}
0 \\ 1
\end{pmatrix} =
\begin{pmatrix}
-1 \\ 1
\end{pmatrix},
\]
the boundary slopes of $N_1\cup N_2\cup N_3\cup N_4\cup_A N_5$ are 
$(0,-1)$. Similarly, after Legendrian surgery along $D^1_{1/2}$ we get 
$N_1\cup N_2\cup N_3\cup_B N_4\cup_A N_5$ with boundary slopes $(0,0)$, 
and after Legendrian surgery along $D^2_{1/2}$ and $D^3_{1/2}$ we get, 
respectively, $N_1\cup N_2\cup N_3\cup_{B^2} N_4\cup_A N_5$ and 
$N_1\cup N_2\cup N_3\cup_{B^3} N_4\cup_A N_5$ with boundary slopes $(0,1)$ 
and $(0,2)$. It is easily checked that 
\[
N_3\cup_{B^3} N_4\cup_A N_5 = 
N_3\cup N_4\cup N'_5,
\]
where $N'_5$ has boundary slopes $(\infty,2)$. After Legendrian surgery 
along $C_{1/4}$ we get 
\[
N_1\cup N_2\cup_A N_3\cup N_4\cup N'_5
= N_1\cup N_2\cup N_3\cup N'_4,
\]
where $N'_4$ has boundary slopes $(0,-2)$. Arguing in a similar
fashion, it is easy to check that after performing the remaining
Legendrian surgeries we end up with $N_1$ glued, in view of
Equation~\eqref{e:relation}, via the identity map to the original
$N_6$. Since $N_2\cup N_3\cup N_4\cup N_5\cong\T_1$, this concludes
the proof.
\end{proof}

\end{document}